%% file: main.tex
\documentclass{IEEEtran4PSCC}
\ifCLASSINFOpdf
   \usepackage[pdftex]{graphicx}
\else
   \usepackage[dvips]{graphicx}
\fi
\usepackage[cmex10]{amsmath}
\ifCLASSOPTIONcompsoc
 \usepackage[caption=false,font=normalsize,labelfont=sf,textfont=sf]{subfig}
\else
 \usepackage[caption=false,font=footnotesize]{subfig}
\fi
\makeatletter
\let\old@ps@headings\ps@headings
\let\old@ps@IEEEtitlepagestyle\ps@IEEEtitlepagestyle
\def\psccfooter#1{%
    \def\ps@headings{%
        \old@ps@headings%
        \def\@oddfoot{\strut\hfill#1\hfill\strut}%
        \def\@evenfoot{\strut\hfill#1\hfill\strut}%
    }%
    \def\ps@IEEEtitlepagestyle{%
        \old@ps@IEEEtitlepagestyle%
        \def\@oddfoot{\strut\hfill#1\hfill\strut}%
        \def\@evenfoot{\strut\hfill#1\hfill\strut}%
    }%
    \ps@headings%
}
\makeatother

\usepackage{microtype}
\usepackage{soul}
\usepackage{xcolor}
\usepackage{tabularx}
\usepackage{multirow}
\usepackage{graphicx}
\usepackage{booktabs}
\usepackage{makecell}

\input{macros}

\begin{document}
%
\title{Refining bridge-block decompositions through two-stage and recursive tree partitioning}

\author{
\IEEEauthorblockN{Leon Lan, Alessandro Zocca}
\IEEEauthorblockA{Department of Mathematics, Vrije Universiteit Amsterdam, NL\\
\{l.lan, a.zocca\}@vu.nl}
}
\maketitle
\begin{abstract}
In transmission networks power flows and network topology are deeply intertwined due to power flow physics. Recent literature shows that specific network substructures named \textit{bridge-blocks} prevent line failures from propagating globally. A two-stage and recursive \textit{tree partitioning} approach have been proposed to create more bridge-blocks in transmission networks, improving their robustness against cascading line failures. In this paper we consider the problem of refining the bridge-block decomposition of a given power network with minimal impact on the maximum congestion. We propose two new solution methods, depending on the preferred power flow model. More specifically, (i) we introduce a novel MILP-based approach that uses the DC approximation to solve more efficiently the second-stage optimization problem of the two-stage approach and (ii) we show how the existing recursive approach can be extended to work with AC power flows, drastically improving the running times when compared to the pre-existing AC-based two-stage method.
\end{abstract}
\vspace{-0.5cm}

\section{Introduction}
Historical blackouts have shown that transmission line failures play an important role in the initiation of cascading failures \cite{Bienstock2015}. The complex network structure of power grids in combination with the underlying power flow physics gives rise to complicated cascading failure patterns, which often exhibit non-local propagation of line failures. A recent paper~\cite{Zocca} shows that specific graph structures, called \textit{bridge-blocks}, ensure that line failures propagate only locally. In particular, it is shown that line failures cannot propagate across lines that act as bridges in the network. These failure localization results have been obtained through studying analytical properties of the line outage distribution factor using the DC approximation, while \cite{BialekJanusz2021,Guo2020a} report similar results when using AC power flows.

\begin{figure}[!ht]
    \centering
    \includegraphics[width=0.95\linewidth,keepaspectratio]{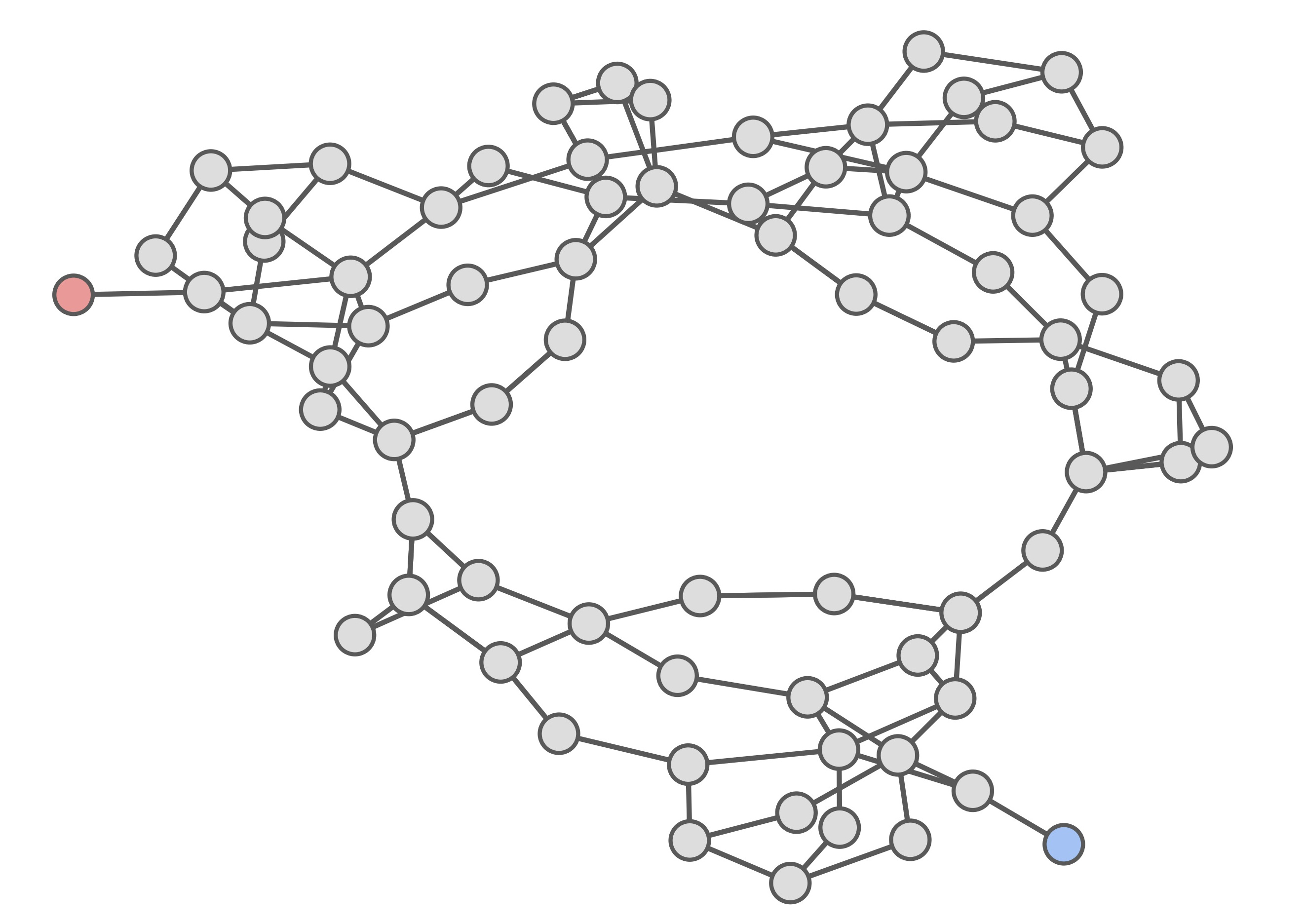}
    \caption{The bridge-block decomposition of the IEEE-73 network, which has one large bridge-block (in gray) of size 71 and two trivial bridge-blocks of size 1 (in other colors). A line failure in the large gray bridge-block could potentially affect any other line that belongs to it.
    }
    \label{fig:IEEE-118_BBD}
\vspace{-0.5cm}
\end{figure}

Despite the potential that bridges may offer in preventing non-local failure propagation, most power grid networks have not been designed with this principle in mind. It is illustrated in~\cite{Zocca} that most test power networks have a \textit{trivial bridge-block decomposition}, i.e., they have one large bridge-block comprising a very large fraction of the network, thus potentially allowing line failures to propagate through the entire network (see Figure~\ref{fig:IEEE-118_BBD}).

Aiming to improve the bridge-block decomposition of a power grid and, hence, its robustness against line failures, \cite{Zocca} proposes to temporarily switch off a set of carefully selected lines, in an adaptive fashion with respect to the current generation and demand patterns. Line switching is a consolidated electricity grid management paradigm that received a lot of attention in the literature. The most well-known applications are \textit{Optimal Transmission Switching} \cite{Salkuti2018,Hedman2011}, where switching actions aim at maximizing economic efficiency of generation dispatch, and \textit{Controlled Islanding} \cite{Ding2013,Quiros-Tortos2015}, where lines are switched off to split the network into disconnected islands as a last-resort emergency measure to stop cascading failures.

A heuristic two-stage procedure called \textit{tree partitioning} has been proposed by \cite{Zocca} and \cite{BialekJanusz2021} (using DC and AC power flow models, respectively) to judiciously switch off lines to improve the bridge-block decomposition of a given power network, while minimizing the impact on line congestion. This two-stage heuristic first identifies suitable clusters and then, by means of switching actions, makes sure that they are connected in a tree-like manner, hence becoming bridge-blocks.
The computational bottleneck of this approach resides in the second stage, since it requires evaluating an exponential number of spanning trees to find the best way to connect the identified clusters. Since the brute-force method does not scale well with large networks, \cite{Zocca} also proposes an alternative recursive approach for tree partitioning, which refines the network bridge-block decomposition at each iteration.

There are two major contributions in this paper. First, we introduce an MILP-based algorithm that solves the second-stage problem \textit{exactly} under the DC power flow model. Our algorithm overcomes the computational bottleneck formed by the existing brute-force algorithms, while producing better results than the recursive algorithm in the considered test cases. Second, aiming to extend tree partitioning to the more realistic AC power flow setting, we modify the recursive approach proposed in~\cite{Zocca} to work with AC power flows. For several test cases, our AC variant of the recursive approach yields solutions qualitatively comparable to those returned by the two-stage brute-force approach under AC power flows, but runs drastically faster and has the potential to scale well for large networks.

The paper is organized as follows. After some preliminaries in Section~\ref{sec:preliminaries}, Section~\ref{sec:problem} formally introduces the optimization problem for bridge block decomposition refinement and revisits two existing approaches to tree partitioning. In Section~\ref{sec:milp}, we outline our MILP-based algorithm for the second stage problem using DC power flows and in Section~\ref{sec:rc} we propose an AC variant of the recursive algorithm. Section~\ref{sec:results} presents the numerical results and we conclude the paper with Section~\ref{sec:conclusion}.

\section{Preliminaries}
\label{sec:preliminaries}
\subsection{Power network model}
We model a transmission network as a connected, directed graph $G=(V,E)$, where $V$ is the set of vertices (buses) and $E$ the set of edges (transmission lines). We denote by \(n=|V|\) the number of buses and by \(m=|E|\) the number of lines. Each bus \(i \in V\) has a \textit{net power injection} \(p_i\), where \(p_i > 0\) is interpreted as \emph{injected} power and \(p_i \leq 0\) as \emph{consumed} power at bus $i$. Each line $\ell=(i,j) \in E$ has a capacity \(c_\ell = c_{ij} > 0\), denoting its rating, i.e., the maximum power that the line can carry.

Throughout this paper, we will use both DC and AC power flow models.
For a full description of AC power flow models, we refer the reader to \cite{Bienstock2015}. We describe here a lossless DC power flow model in which generation always matches demand, i.e., $\sum_{i=1}^{n} p_i = 0$. We refer to any such vector $\bm{p}$ of power injections as \textit{balanced}. Let $\theta_{i} \in \bbR$ denote the \textit{phase angle} of bus $i$. For each line $\ell = (i,j)$, let $f_\ell = f_{ij} \in \bbR$ denote the active power flow and let \(b_\ell = b_{ij} > 0\) denote the \textit{line susceptance}. Given a vector of power injections $\bm{p} \in \bbR^{n}$, the corresponding line flows $\bm{f} \in \bbR^{m}$ and phase angles $\bm{\theta} \in \bbR^{n}$ are obtained by solving the \textit{DC power flow equations}:
\begin{subequations}
\label{eq:DCPF}
\begin{align}
p_i &= \sum_{j: (i, j) \in E} f_{ij} - \sum_{j: (j, i) \in E} f_{ji}, & \forall \, i \in V, \label{eq:DCPF-1}\\
  f_{ij} &=  b_{ij}(\theta_i - \theta_j), & \forall \, (i, j) \in E. \label{eq:DCPF-2}
\end{align}
\end{subequations}
Equation~\eqref{eq:DCPF-1} guarantees flow conservation and \eqref{eq:DCPF-2} captures the flow dependency on susceptance and angle differences. The DC power flow equations \eqref{eq:DCPF} admit a unique power flow solution $\bm{f}$ for each balanced injection vector \(\bm{p}\). 
The solution for the phase angles $\bm{\theta}$ is unique up to an arbitrary reference angle: without loss of generality, we select bus $n$ as the \textit{reference bus} with phase angle $\theta_n = 0$. 

\subsection{Tree partitions and bridge-block decomposition}
We now briefly review some graph-theoretical terminology.
A \textit{$k$-partition} of the network $G$ is a collection $\calP = \{ \calV_1, \calV_2, \dots, \calV_k \}$ of non-empty, disjoint vertex sets $\calV_1, \calV_2, \dots, \calV_k$ called \textit{clusters} such that $\bigcup_{i=1}^k \calV_k = V$.
We denote by $\nodeset{k} = \{1, 2, \dots, k\}$ the set of integers from $1$ to $k$, which will be used to denote the clusters. Given a partition $\calP$, a line $(i,j)$ is called an \textit{internal edge} if both $i$ and $j$ belong to the same cluster and a \textit{cross edge} otherwise. We denote by $E_C(\calP)$ the set of cross edges determined by partition $\calP$. The \textit{reduced graph} $G_\calP$ is the graph whose vertices are the clusters in $\calP$ and where an edge is drawn for each cross edge connecting two different clusters. Note that it is possible for the reduced graph to have multiple edges between two vertices, and thus to be a multigraph.

\begin{figure}[!b]
    \centering
    \subfloat[]{{\includegraphics[width=0.48\linewidth]{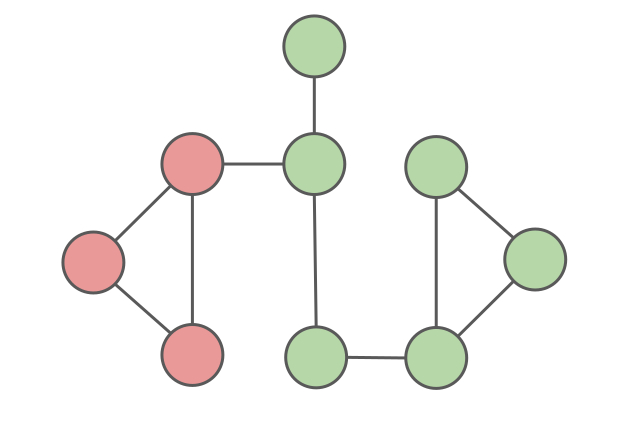} \includegraphics[width=0.48\linewidth]{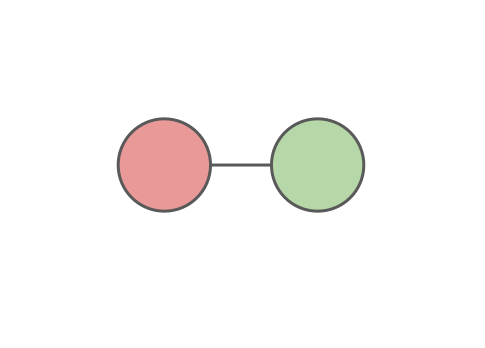}      \label{fig:tree-partition}}}%
    \vspace{0.5cm}
    \subfloat[]{{\includegraphics[width=0.48\linewidth]{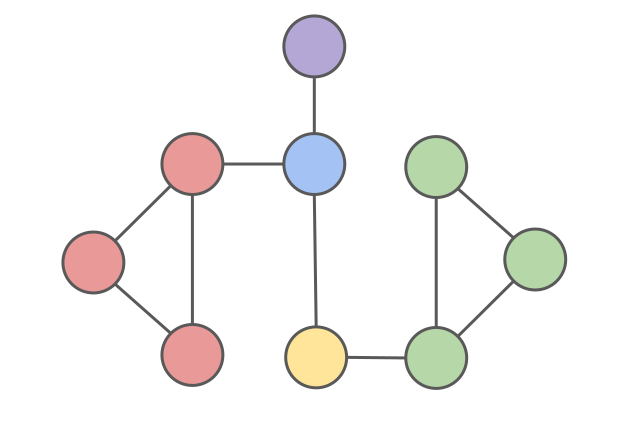}
    \includegraphics[width=0.48\linewidth]{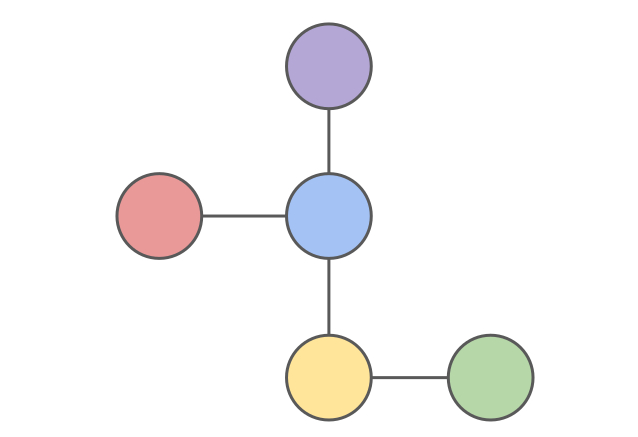}
    \label{fig:bbd}}}
    \vspace{0.5cm}
    \caption{(a) A tree partition of a graph $G$ and the corresponding reduced graph. (b) The bridge-block decomposition of $G$ and the corresponding reduced graph. Note that, being the bridge-block decomposition, it is at least as fine as any tree partition of the same network and in fact is finer than (a).
    }%
    \label{fig:tp-bbd}
\end{figure}

We say that a partition $\calP$ is a \textit{tree partition} of $G$ if the reduced graph $G_\calP$ is a tree (see Figure~\ref{fig:tree-partition}). A \textit{bridge} is a cut-edge for the graph i.e., an edge whose removal would disconnect it. The \textit{bridge-block decomposition} of a graph is its partition in disconnected subgraphs that is obtained after removing all bridges in the graph (see Figure~\ref{fig:bbd}). Each cluster in the bridge-block decomposition is called a \textit{bridge-block}. It is easy to show that the bridge-block decomposition is a tree partition, and, in particular, is the \textit{irreducible} one, i.e., the maximal by inclusion~\cite{Guo2018}. Consequently, there could exist tree partitions that are not the bridge-block decomposition, as shown by comparing Figure~\ref{fig:tree-partition}~and~\ref{fig:bbd}.
Given any tree partition $\calP$, we henceforth say that the bridge-block decomposition is always \textit{as fine} as $\calP$, meaning that every bridge-block is contained in some cluster of the tree partition. 

\begin{figure*}[!t]
    \centering
    \subfloat[]{\includegraphics[width=0.31\linewidth]{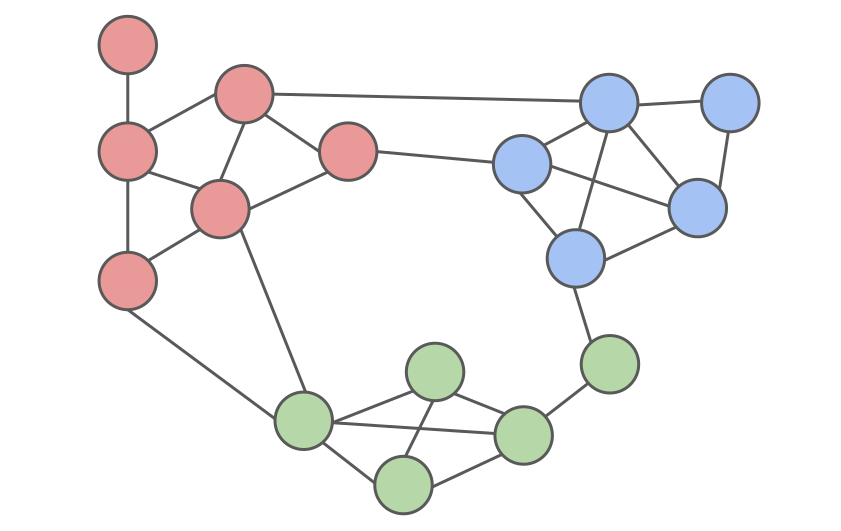}%
    \label{fig_first_case}}
    \hfill
    \subfloat[]{\includegraphics[width=0.31\linewidth]{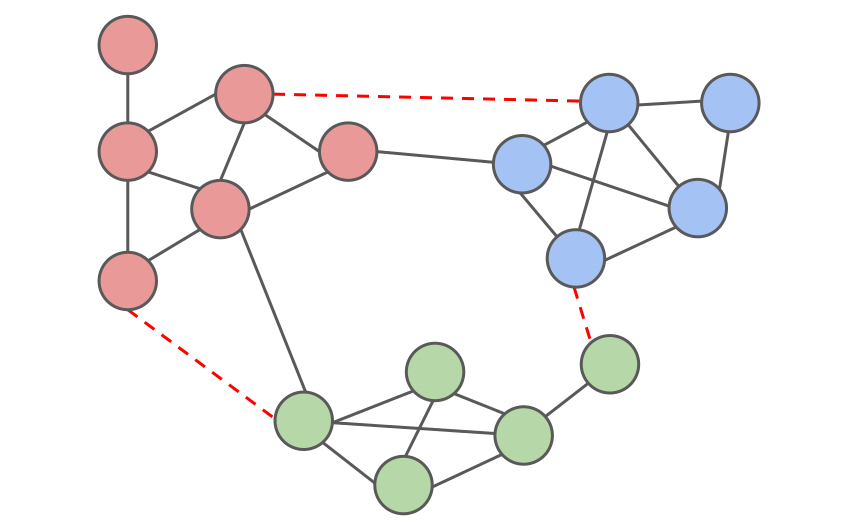}%
    \label{fig_second_case}}
    \hfill
    \subfloat[]{\includegraphics[width=0.31\linewidth]{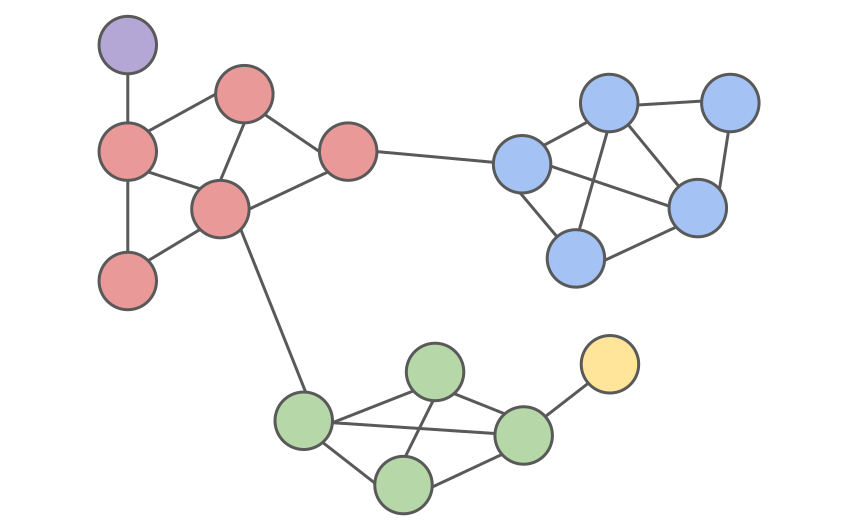}
    \label{fig_third_case}}
    \caption{Illustration of the two-stage approach for $k=3$ target clusters. (a) A 3-partition $\calP$ of a network $G$. (b) A subset of lines $\calE$ (dashed, in red) is switched off, turning $\calP$ into a tree partition of $G^\calE$. (c) The resulting bridge-block decomposition of $G^\calE$, which is slightly finer than the identified partition $\calP$.
}
    \label{fig:two-stage}
    \vspace{-0.45cm}
\end{figure*}

\section{Bridge-block decomposition refinement}
\label{sec:problem}
\subsection{Motivation}
After the failure or disconnection of a transmission line, its original power flow gets \textit{globally} redistributed as prescribed by power flow physics on the remaining lines, some of which can overload and also get disconnected. We refer to this phenomenon as \textit{failure propagation}.
It was shown in \cite{Zocca,BialekJanusz2021,Guo2018,Guo2020} that a line failure does not propagate across bridges, but instead only impacts the flows on lines that belong to the same bridge-block. Most power networks, however, have a very meshed structure and trivial bridge-block decompositions~\cite{Zocca}, making them very prone to non-local line failure propagation, see, e.g., Figure~\ref{fig:IEEE-118_BBD}.

The number of (non-trivial) bridge-blocks can be increased by switching off lines in a procedure named \textit{bridge-block decomposition refinement} \cite{Zocca}. Nevertheless, as there are an exponential number of lines to consider, there is no obvious way to select which lines to remove in general. \cite{Zocca} proposes a bottom-up approach, named \textit{tree partitioning}: a target partition is first identified using clustering methods, and then these clusters are ensured to be connected in a tree-like manner, transforming the identified partition into a tree partition. More formally, given a power network $G=(V,E)$, the goal of tree partitioning is to identify a partition $\calP$ and a subset of lines $\calE$ to be switched off, such that $\calP$ is a tree partition of the post-switching network $G^\calE = (V, E \setminus \calE)$.

This tree partitioning procedure guarantees that the post-switching network has a bridge-block decomposition that is at least as fine as the tree partition $\calP$.

\subsection{Line congestion}
Let $g_\ell$ denote the \textit{congestion level} on line $\ell$, whose specific formula depends on the used power flow model. Under the DC approximation we define the congestion level as the non-negative ratio $g_\ell:=\abs{f_\ell} \mathbin{/} c_{\ell}$.
Under the AC power flow model the calculation of the congestion is more involved and the apparent power should be considered \cite{BialekJanusz2021}. In either case, we say that a line $\ell$ is congested if $g_\ell > 1$.

It is desirable that the switching actions do not cause any of the remaining lines to become congested as a result of flow redistribution, but this is hard to ensure up front due to the complexity of the power network and the underlying power flow physics. We thus select the set $\calE$ of lines to be switched off that minimizes the \emph{maximum congestion}, defined as
\begin{align}
  \gamma(\calE) \defeq \max_{\ell \in E\setminus \calE} \widetilde{g}_\ell,
\end{align}
where $\widetilde{g}_\ell$ is the post-switching congestion level on line $\ell$, calculated using the power flow equations on $G^\calE$ assuming that the power injections $\bm{p}$ are unchanged. Note that keeping the maximum congestion below 1 directly implies that the post-switching network has no congested lines. Furthermore, line flows may slightly exceed the capacities as long as subsequent remedial actions are undertaken to alleviate the congestion \cite{BialekJanusz2021}.%
\subsection{Problem statement}
We now formulate an optimization problem to tree partition a network with minimal impact on the maximum congestion based on \cite{Zocca}~and~\cite{BialekJanusz2021}.
Given a power network $G=(V,E)$ with balanced power injections $\bm{p}$ and a positive integer $k \geq 2$, the goal is to identify a $k$-partition $\calP$ and a set of lines $\calE \subset E$ to be switched off such that we minimize the maximum congestion $\gamma(\calE)$ and satisfy the following properties:
\begin{itemize}
    \item the post-switching network $\GcalE=(V, E \setminus \mathcal{E})$ is connected;
    \item $\calP$ is a tree partition of $G^\calE$.
\end{itemize}

Determining the optimal number $k$ of clusters is outside the scope of this paper, and we henceforth assume that $k$ is a given input parameter. Note further that the problem formulation does not specify any minimum size for the identified clusters in $\calP$. However, in view of the considered spectral methods that intrinsically use normalization, obtaining a partition with trivial clusters is a rather rare occurrence.

\begin{figure*}[!t]
    \centering
    \subfloat[\label{fig:rc1-1}]{\includegraphics[width=0.31\linewidth]{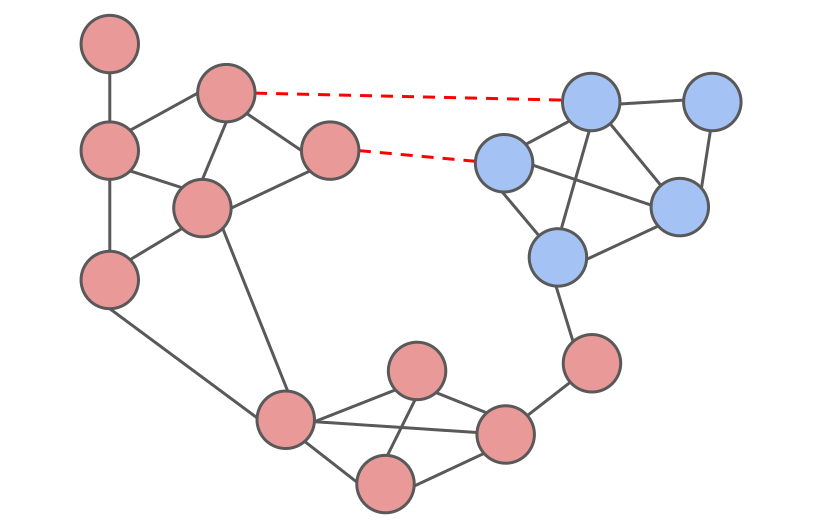}%
    }
    \hfill
    \subfloat[\label{fig:rc1-2}]{\includegraphics[width=0.31\linewidth]{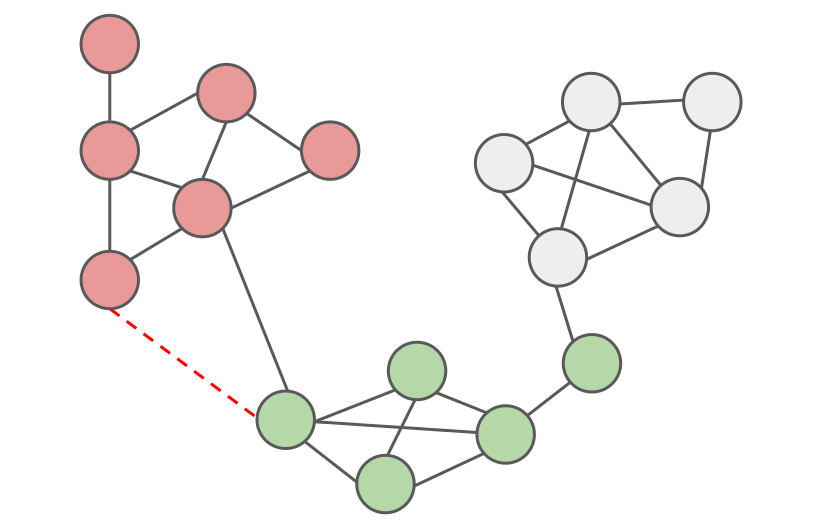}%
    }
    \hfill
    \subfloat[\label{fig:rc1-3}]{\includegraphics[width=0.31\linewidth]{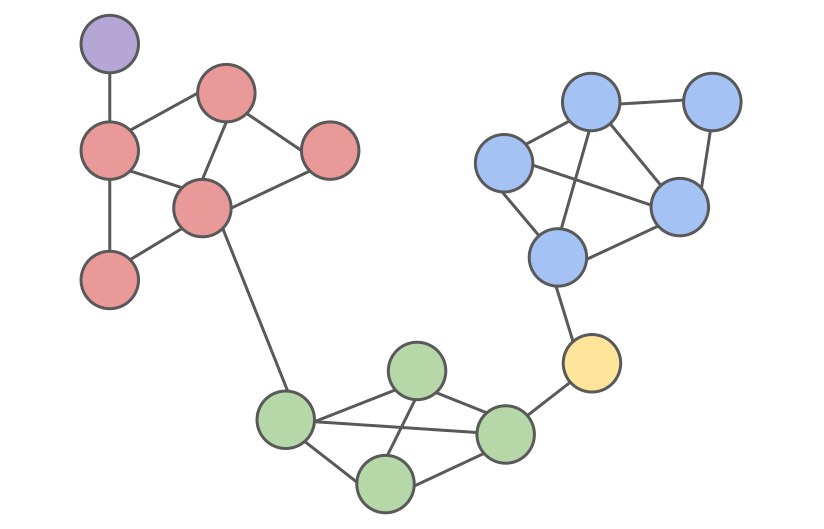}
    }
    \caption{Illustration of the recursive approach for $k=3$ target clusters. (a) The first iteration refines the network into a large red bridge-block and a smaller blue bridge-block. (b) The second iteration refines the largest bridge-block from the previous iteration, i.e., the red one in (a). (c) The resulting post-switching network after $k-1=2$ iterations and its bridge-block decomposition. Note that the latter is identical to the post-switching bridge-block decomposition obtained using the two-stage approach (cf.~Figure~\ref{fig:two-stage}), yet the post-switching network and hence the congestion levels are slightly different.
}
    \label{fig:rc1}
    \vspace{-0.45cm}
\end{figure*}

\subsection{Two-stage approach}
\label{sec:two-stage}
We describe the two-stage heuristic approach proposed by~\cite{Zocca} to tree partition a network with minimal impact on the congestion. These two stages arise naturally because identifying $\calP$ and $\calE$ simultaneously is computationally intractable for the optimization problem at hand. The two-stage approach decouples them: the first stage computes the partition $\calP$ and the second stage returns the corresponding best subset $\calE$. 

\subsubsection{First stage -- Optimal Bridge-blocks Identification}%
\label{sec:obi}%
For a given $k$, the first stage aims to find a sensible partition of the power network into $k$ clusters so that the identified partition can then be transformed into a tree partition in the second stage. A ``good'' partition must have few cross edges between clusters and with very modest power flows. Indeed, switching off too many lines or lines with large power flows often leads to large power flow redistribution and results in severe congestion in the resulting network. We thus seek a partition that has (i) very few (and/or low weight) cross edges between the clusters and (ii) clusters that are similar in size and balanced in terms of total net power. In~\cite{Zocca}, the authors formulate this as an optimization problem, named the \textit{Optimal Bridge-blocks Identification (OBI)} problem, using a power-flow-weighted version of the \emph{network modularity} problem \cite{Newman2006} restricted to $k$-partitions.

The modularity maximization problem is NP-hard \cite{Brandes2006} and thus so is the OBI problem, hence two strategies are considered in~\cite{Zocca} to obtain good approximate solutions.
The first one is spectral clustering, using either the normalized Laplacian matrix $L_N$ or the normalized modularity matrix $B_N$ of the network with absolute power flows as edge weights.
Both these variants approximate the solution of the same problem: \cite{Yu2010} shows that the optimal solution to the normalized cut problem for a fixed number of target clusters $k$ (i.e., spectral clustering using $L_N$) is identical to that of the normalized modularity problem (i.e., spectral clustering using $B_N$). However, due to differences in their eigensystems, both algorithms often yield slightly different results and thus we consider them both separately as \textit{Spectral $L_N$} and \textit{Spectral $B_N$} for our numerical results.
The second strategy employs the \textit{Fastgreedy} algorithm \cite{Clauset2004}, which is a fast heuristic algorithm that creates clusters in a hierarchical agglomerate way.

\subsubsection{Second stage -- Optimal Bridge Selection}
In \cite{Zocca}, the second stage is formulated as an optimization problem, named \emph{Optimal Bridge Selection (OBS)}, where the identified partition from the first stage is turned into a tree partition on the post-switching network. Consider the power network $G=(V,E)$ and the $k$-partition $\calP$ obtained by solving the OBI problem. The goal of the OBS problem is to remove a subset of cross edges $\calE \subset E_C(\calP)$ such that we minimize the maximum congestion $\gamma(\calE)$ on the post-switching network $G^\calE=(V, E \setminus \calE)$, where we assume that the power injections $\bm{p}$ are unchanged, and such that $\calP$ is a tree partition of $G^\calE$. 

An alternative formulation of the OBS problem, which will be leveraged for the MILP-based algorithm, can be given using an explicit definition of the reduced graph. We define the reduced graph of $G$ corresponding to $\calP$ as $G_{\calP}=(\nodeset{k}, E_{C}(\calP))$, i.e., as the graph whose vertices are the clusters indexed from $1$ to $k$ and whose edges are the cross edges between them. Solving the OBS problem then corresponds to finding a spanning tree $T$ on the reduced graph $G_\calP$, such that the removal of lines $\calE = E_C(\calP) \setminus T$ minimizes the maximum congestion $\gamma(\calE)$ on the post-switching network $G^\calE$.

The OBS problem is particularly difficult since, for any given subset $\calE$, we need to recalculate the power flows on the post-switching graph $G^\calE$ in order to obtain the the maximum congestion $\gamma(\calE)$. In previous studies \cite{Zocca,BialekJanusz2021}, the OBS problem has been solved using brute-force algorithms that enumerate all spanning trees. However, since their number is exponential in the number of clusters $k$, any brute force approach is intractable for large instances.

\subsection{Recursive approach}
\label{sub:ra}

The hardness of the OBS problem limits the applicability of the two-stage approach to larger network instances.
However, if the selected partition $\calP$ consists of only $k=2$ clusters, the number of spanning trees in the reduced graph $G_\calP$ is exactly equal to the number of cross edges. Therefore, the OBS problem can be solved much faster in this case than when considering $k>2$ clusters at once. This key idea is at the core of the \emph{recursive approach} to tree partitioning, introduced by \cite{Zocca} precisely to overcome the hardness of the OBS problem.

In the recursive approach, the largest bridge-block in the network is iteratively refined into two smaller bridge-blocks until the desired number of clusters is obtained. In other words, at every iteration of the recursive algorithm, one first solves the OBI problem restricted to $2$-partitions and then solves the OBS problem given the resulting bipartition. Hence, given the number $k$, the recursive algorithm solves the OBI and OBS problem for a total of $k-1$ times. Figure~\ref{fig:rc1} demonstrates how the recursive approach works for a small network.

The main advantage of the recursive approach is that is much faster than the two-stage approach, having time complexity linear in $k$. 
The OBS problem can now be solved using a brute-force algorithm that considers a linear number of spanning trees on the reduced graph and computes the corresponding power flows and congestion levels on the post-switching network. Consequently, the running time of the recursive approach mainly depends on the speed of the clustering algorithm and the number of clusters $k$ to be obtained.

Even when starting with the same network, the solutions from the two-stage and recursive approaches could differ from each other as shown by comparing Figure~\ref{fig:two-stage}~and~\ref{fig:rc1}, and hence lead to different congestion levels. 
Later in Section~\ref{sec:results} we will compare the performance between these two approaches (both using DC and AC power flows).

\section{MILP formulation for the OBS problem using DC power flows}
\label{sec:milp}
In this section, we present our first main contribution where we show how to solve the OBS problem using an exact MILP-based algorithm that uses the DC power flow model. The MILP-based algorithm overcomes the bottleneck imposed by the earlier proposed brute-force algorithms, consequently allowing the two-stage approach to be solved more efficiently.

\begin{subequations}
\label{opt:obs}
Assuming the partition $\calP$ is given, we start by describing in detail how the OBS problem can be cast into an MILP.
We say that a line is \emph{inactive} if it is switched off and \emph{active} otherwise. Let the decision variable $\gamma \in \bbR$ denote the maximum congestion and let the decision variables $f_{ij} \in \bbR, (i, j) \in E$ denote the active power flow on the lines. The objective is to minimize the maximum congestion in the network:
\begin{align}
    \min \quad \gamma,
\end{align}
which is subject to
\begin{align}
\gamma \geq \abs{f_{ij}} / c_{ij}, &\qquad \forall (i, j) \in E.
\end{align}

First, we introduce constraints to obtain a spanning tree on the reduced graph. With slight abuse of notation, we henceforth identify a cross edge $\ell \in \crossE$ with indices $(i,j)$ when we refer to the line in the original network $G$ and $(u,v)$ when we refer to the reduced graph $G_\calP$, where $i,j$ are bus indices and $u,v$ are cluster indices. Recall that $k$ is the number of clusters in $\calP$ and thus also the number of vertices in the reduced graph. Let the binary decision variables $y_{uv} \in \binaryset, (u, v) \in E_C(\calP)$ indicate whether or not a cross edge is active. We introduce a cardinality constraint on the number of cross edges:
\begin{align}
    \label{eq:obs_st}
    \sum_{(u, v) \in E_C(\calP)} y_{uv} = k-1.
\end{align}
Moreover, to ensure that the post-switching reduced graph is connected, we add the \textit{single commodity flow} constraints \cite{Gavish1978}. The main idea is to assign a source vertex that sends a fictitious unit flow to all other vertices, which is possible if and only if the graph is connected. Let the decision variables $q_{uv} \in \bbR, (u, v) \in E_C(\calP)$ denote the commodity flow on the cross edges and we assign vertex $1$ as the source vertex. Then, the single commodity flow constraints are expressed as follows:
{\fontsize{10}{0}
\begin{align}
    \, \smashoperator{\sum_{(1, v) \in E_C(\calP)}} \, \, q_{1v} \,  - \quad \, \smashoperator{\sum_{(v, 1) \in  E_C(\calP)}} \, \, q_{v1} \, = k-1,\label{eq:obs_scf_1}  & \\
    \smashoperator{\sum_{(u, v) \in E_C(\calP)}} \, \,  q_{uv} \, - \quad \, \smashoperator{\sum_{(v, u) \in E_C(\calP)}} \, \,  q_{vu} \, = -1, & \quad \forall \, u \in \nodeset{k} \setminus \{1\},\label{eq:obs_scf_2}  \\
    -(k-1) y_{uv} \leq q_{uv} \leq (k-1) y_{uv}, & \quad \forall \, (u, v) \in E_C(\calP). \label{eq:obs_scf_3}
\end{align}}%
Equation~\Cref{eq:obs_scf_1} ensures that the net commodity flow, defined as the outgoing minus the incoming commodity flow, of the source vertex is exactly $k-1$, i.e., it produces $k-1$ units of commodity flow. Similarly, \cref{eq:obs_scf_2} ensures that the net commodity flow of the demand vertices is $-1$, meaning that they consume one unit of commodity flow. Equation~\cref{eq:obs_scf_3} ensures that the inactive lines carry no commodity flow. Hence, since \cref{eq:obs_st,eq:obs_scf_1,eq:obs_scf_2,eq:obs_scf_3} ensure that there are $k-1$ cross edges and the reduced graph is connected, this implies that we obtain a spanning tree on the reduced graph.

Next, we model the impact of the switching actions in terms of congestion on the post-switching network. Let the decision variables $\theta_i \in \bbR, i \in V$ denote the phase angle of each bus. We first consider the cross edges:
{\fontsize{10}{0}
\begin{align}
f_{ij} \leq b_{ij}(\theta_i - \theta_j) + M_{uv}(1-y_{uv}), &\quad \forall \, \ell \in \crossE, \label{eq:dcopf_switch_1}\\
f_{ij} \geq b_{ij}(\theta_i - \theta_j) + M_{uv}(1-y_{uv}), &\quad \forall \, \ell \in \crossE, \label{eq:dcopf_switch_2}\\
f_{ij} \leq M_{uv}y_{uv}, & \quad \forall \, \ell \in \crossE, \label{eq:dcopf_switch_3}\\
f_{ij} \geq M_{uv}y_{uv}, & \quad \forall \, \ell \in \crossE. \label{eq:dcopf_switch_4}
\end{align}}%
Equations~\cref{eq:dcopf_switch_1,eq:dcopf_switch_2} ensure that the DC power flow equations hold on the active cross edges, whereas \cref{eq:dcopf_switch_3,eq:dcopf_switch_4} switch off the DC power flow equations for all inactive cross edges. We set the big-M value $M_{uv}$ at four times the capacity of the corresponding line. Moreover, we add the following constraints:
\begin{align}
f_{ij} = b_{ij}(\theta_i - \theta_j), & \quad \forall \, (i,j) \in E \setminus \crossE,\label{eq:dcopf_internal_edges}\\
\sum_{(i,j) \in E} f_{ij} - \sum_{(j, i) \in E} f_{ji} = p_i, &\quad \forall \, i \in V, \label{eq:flow_conservation}\\
\theta_n = 0. \label{eq:refbus}
\end{align}
Equation~\eqref{eq:dcopf_internal_edges} models the DC power flows equations on the internal edges, which are all active by assumption, and \cref{eq:flow_conservation,eq:refbus} ensure flow conservation and assign $n$ as the reference bus.
\end{subequations}

In summary, the MILP formulation \eqref{opt:obs} contains (i) a spanning tree formulation on the reduced graph and (ii) DC power flow equations to compute the congestion on the post-switching network. Assuming the partition is obtained solving the OBI problem (cf.~Section~\ref{sec:two-stage}), the number of cross edges is relatively small, meaning that modeling the spanning tree formulation requires a relatively small number of constraints and decision variables. Hence, the size of the MILP depends mostly on the network instance size $|V|$ and $|E|$. In Section~\ref{sec:dc_results}, we compare the performance between the newly improved two-stage approach and recursive approach assuming DC power flows.

\section{AC modification of the recursive approach}
\label{sec:rc}
The objective function of the OBS problem requires to calculate the congestion $\gamma(\calE)$ in the post-switching network $G^\calE$ for all possible sets $\calE$ of switching actions. Ideally one would calculate line flows and thus the maximum congestion using an AC power flow model, but this is not feasible when solving the OBS problem with brute force~\cite{BialekJanusz2021,Zocca} or with the MILP approach proposed in the previous section. Both these methods strongly rely on the linear DC power flow model.
The recursive approach, which was originally introduced to solve the tree partitioning problem using the DC power flow model, does not suffer from the same scalability issue and it is thus possible to incorporate AC power flow calculations as subroutine.
This simple yet crucial modification enables the recursive approach to solve the tree partitioning problem with AC power flows, being the first algorithm to do so with time complexity that is linear in $k$. In Section~\ref{sec:ac_results}, we compare the performance between the recursive and two-stage approach using AC power flows, where we use a brute-force algorithm to solve the OBS problem in the two-stage approach.

\section{Numerical results}
\label{sec:results}
In this section, we evaluate the performance of the discussed methods to solve the optimization problem introduced in Section~\ref{sec:problem}.
In our numerical experiments we use test cases from the \texttt{PGLib-OPF} library \cite{Babaeinejadsarookolaee2021}.
The experiments are performed using an Intel® Core™ i7-8750H CPU @ 2.20GHz × 12 and 16 GB RAM.
For more implementation details see \cite{Lan}.

\begin{table}[!b]
\centering
\caption{Bridge-block decomposition due to tree partitioning with $k=5$ target bridge-blocks.}
\label{tab:bbdr}
\begin{tabular}{@{}lcrlcr@{}}
\toprule
          & \multicolumn{2}{c}{\thead{pre-switching\\ non-trivial \\bridge-blocks}}&  & \multicolumn{2}{c}{\thead{post-switching\\ non-trivial\\ bridge-blocks}}\\
          \cmidrule{2-3} \cmidrule{5-6}
Case      & \thead{\#} & \multicolumn{1}{r}{sizes} &  & \thead{\#} & \multicolumn{1}{r}{sizes largest $5$}    \\ \midrule
IEEE-30   & 1                                                & $\{27\}$                       &  & 5                                                & \{$7, 3, 3, 3, 3\}$           \\
IEEE-118  & 1                                                & $\{109\}$                      &  & 5                                                & \{$39, 19, 19, 9, 8\}$        \\
GOC-179   & 1                                                & $\{150\}$                      &  & 8                                                & \{$40, 30, 30, 15, 13\}$      \\
ACTIV-200 & 3                                                & $\{125, 4 ,3\}$                &  & 6                                                & \{$37, 20, 15, 12, 4\}$       \\
IEEE-300  & 3                                                & $\{206, 3, 3\}$                &  & 8                                                & \{$58, 46, 36, 19, 16\}$      \\
GOC-500   & 3                                                & $\{354\}$                      &  & 5                                                & \{$92, 84, 59, 53, 28\}$      \\
GOC-793   & 1                                                & $\{500\}$                      &  & 6                                                & \{$96, 61, 53, 47, 41\}$      \\
RTE-1888  & 2                                                & $\{881, 5\}$                   &  & 8                                                & \{$228, 165, 158, 146, 131\}$ \\
\bottomrule
\end{tabular}
\end{table}
\subsection{Post-switching bridge-block decompositions}
\label{sec:bbd_results}
In this section we quickly look at the quality of the bridge-block decompositions obtained by tree partitioning various test networks. We use the two-stage approach on several test networks with $k=5$ target clusters, where the OBI problem is solved using $\SpectralLn$ and the OBS problem is solved by the MILP-based algorithm. Table~\ref{tab:bbdr} shows the bridge-block decomposition characteristics of the pre-switching and post-switching networks. In particular, observe that the bridge-block decomposition of the pre-switching networks all consist of a single bridge-block encompassing a large fraction of the buses. The tree partitioning procedure creates post-switching networks that have \textit{at least} 5 non-trivial bridge blocks, but often slightly more due to unintended formation of new bridges (e.g., see Figure~\ref{fig:rc1}). The new non-trivial bridge-blocks are smaller in size, and, consequently, the new networks are more robust against failure propagation.
Similar results with respect to the post-switching bridge-block decomposition have been obtained when using other clustering algorithms, i.e., Fastgreedy and $\SpectralBn$, as well as when using the recursive method.

\begin{table}[!hb]
\centering
\caption{Comparison between the MILP-based two-stage and recursive approach with $k=5$ and using DC power flows.}
\label{tab:rc_vs_milp_k=5}
\begin{tabular}{llcccHHccc}
\toprule
 & & \multicolumn{2}{c}{$\gamma(\calE)$} & & & & &\multicolumn{2}{c}{Running time (s)} \\
\cmidrule{3-4} \cmidrule{9-10}
Case & $ALG(\calP)$ & MILP & R-DC & & R-DC & MILP & & MILP & R-DC \\
\midrule
IEEE-118  & Fastgreedy  & 1.57     & \textbf{1.14}      &  & 4     & 7       &  & 0.23   & 0.49  \\
          & $\SpectralBn$  & 1.78  & \textbf{1.00}        &  & 2     & 11      &  & 0.31 & 1.13    \\
          & $\SpectralLn$   & \textbf{1.21}  & \textbf{1.21}         &  & 3     & 3       &  & 0.26   & 0.69  \\
\midrule
GOC-179   & Fastgreedy  & \textbf{1.38}     & \textbf{1.38}       &  & 6     & 5       &  & 0.51  & 0.19   \\
          & $\SpectralBn$ & \textbf{1.38}     & \textbf{1.38}      &  & 6     & 7       &  & 0.93  & 0.26   \\
          & $\SpectralLn$  & \textbf{1.24}  & 1.51         &  & 6     & 7       &  & 0.70  & 0.35   \\
\midrule
IEEE-300  & Fastgreedy  & \textbf{1.16}    & 1.20        &  & 7     & 5       &    & 0.50 & 0.76  \\
          & $\SpectralBn$    & \textbf{1.09}   & 1.68      &  & 10    & 5     &    & 0.59    & 1.12 \\
          & $\SpectralLn$     &\textbf{1.09}    & 1.22     &  & 7     & 4       &   & 0.42 & 0.81   \\
\midrule
GOC-500   & Fastgreedy     & \textbf{1.28}  & 2.38      &  & 15    & 13      &   & 2.93& 1.26    \\
          & $\SpectralBn$    & \textbf{1.01}    & 2.36     &  & 14    & 19      &  & 1.16   & 1.83  \\
          & $\SpectralLn$     & \textbf{1.01} & 2.39        &  & 11    & 17      &  & 1.27  & 1.07   \\
\midrule
GOC-793   & Fastgreedy      & \textbf{1.50}  & 1.54      &  & 28    & 12      &  & 5.88   & 1.67  \\
          & $\SpectralBn$    & \textbf{1.44} & 2.64        &  & 32    & 19      &  & 3.09    & 1.59 \\
          & $\SpectralLn$      & 1.79   & \textbf{1.34}     &  & 14    & 18      &   & 3.82  & 1.55  \\
\midrule
RTE-1888  & Fastgreedy    & \textbf{1.00}    & 1.88      &  & 10    & 0       &   & 4.56  & 5.20  \\
          & $\SpectralBn$   & \textbf{1.00}   & 1.06       &  & 1     & 0       & & 5.51    & 8.56  \\
          & $\SpectralLn$     & 1.10    & \textbf{0.86}     &  & 0     & 5       &  & 12.81   & 4.36 \\
\bottomrule
\end{tabular}
\vspace{-0.35cm}
\end{table}

\subsection{Two-stage vs. recursive approach: DC power flows}
\label{sec:dc_results}
Having demonstrated how tree partitioning can improve the quality of a network's bridge-block decomposition, we now investigate its impact on the maximum congestion under the assumption of DC power flows. Specifically, we compare the performance of the two-stage and recursive approach, where the former is solved by using the MILP-based algorithm for the OBS problem. We denote the two approaches by MILP and R-DC, respectively. We select $k=5$ to demonstrate the scalability of our MILP-based algorithm and use the three different clustering algorithms as presented in Section~\ref{sec:obi}, which we denote by $ALG(\calP)$.

The algorithms were run on a large collection of test networks. For each of them, the initial power injections $\bm{p}$ and power flows $\bm{f}$ were obtained by solving a DC-OPF problem. All presented test cases are selected such that the pre-switching maximum congestion is exactly $1$, i.e., no line carries more power flow than its capacity.

Table~\ref{tab:rc_vs_milp_k=5} reports the objective value and the running time of both MILP and R-DC. For each test case and partitioning method, we highlighted the best solution in terms of maximum congestion in bold. MILP computed the best solutions in most of the presented test cases, often with much lower congestion levels than R-DC. The running times of both methods were comparable for smaller instances, but for large instances MILP took up to 2-4 times longer than R-DC to compute.

The results show that MILP performs better than 2-RC in terms of maximum congestion. Most importantly, MILP has running times comparable to R-DC, meaning that MILP drastically improves the running times of earlier proposed brute-force algorithms for the OBS problem and allows the two-stage approach to be used for large instances. Lastly, the experiments show there is no single best clustering algorithm to be used for tree partitioning. This is partly due to the fact that the OBI problem does not take into account other important network characteristics, such as pre-switching congestion levels and line susceptances. Future work should look into possibly new formulations that take these factors into account.

\subsection{Two-stage vs. recursive approach: AC power flows}
\label{sec:ac_results}
We now compare the performance of the two-stage and recursive approach when using AC power flows. Since the MILP-based formulation of the OBS problem cannot account for the nonlinear AC power flows, we consider the brute-force (BF) variant of the two-stage approach here instead. We henceforth denote the recursive approach by R-AC.

Table~\ref{tab:ac_results} reports the performance of BF and R-AC in terms of maximum congestion and running time on five test cases, where we again highlighted the best solution in terms of maximum congestion in bold. For each test case, the initial power injections and power flows were computed using AC-OPF, but note that the pre-switching congestion $\gamma(\emptyset)$ was not strictly 1. The two-stage and recursive approach show similar performance in terms of maximum congestion: roughly one half of the best solutions were produced by BF, whereas the other half were computed by R-AC. However, most importantly, the recursive approach runs extremely fast with sub-second running times for all considered cases. Lastly, no single clustering algorithm worked best in minimizing the maximum congestion.

\begin{table}[]
\centering
\setlength{\tabcolsep}{0.61em} 
\caption{Comparison between the brute-force-based two-stage and recursive approach with $k=5$ and using AC power flows.}
\label{tab:ac_results}
\begin{tabular}{@{}lHlcccccc@{}}
\toprule
\multicolumn{4}{l}{} & \multicolumn{2}{c}{$\gamma(\calE)$} & & \multicolumn{2}{c}{Running time (s)} \\
\cmidrule{5-6} \cmidrule{8-9}
 Case        & $k$ \, & $ALG(\calP)$   &$\gamma(\emptyset)$ & BF &  R-AC &&  BF &  R-AC\\
\midrule
IEEE-30   & 4                     & Fastgreedy  & 1.07                            & \textbf{1.02}     & 2.13 &      &  1.53                     & 0.05                             \\
   & 4                     & $\SpectralLn$  & 1.07                            & \textbf{1.02}     & 2.13 &      & 1.65   & 0.06                                                \\
   & 4                     & $\SpectralBn$ & 1.07                            & \textbf{1.02}    & 1.06  &  & 1.84   & 0.08                                            \\\midrule
EPRI-39   & 4                     & Fastgreedy  & 0.89                            & \textbf{1.11}     & \textbf{1.11} &     & 0.86       & 0.05                                             \\
   & 4                     & $\SpectralLn$  & 0.89                            & \textbf{0.82}     & 1.11  &              & 0.50   & 0.07                                       \\
   & 4                     & $\SpectralBn$ & 0.89                            & 1.11    & \textbf{1.09}   &          & 1.34    & 0.07                                          \\\midrule
IEEE-73   & 4                     & Fastgreedy  & 0.95                             & 0.96  & \textbf{0.95}   &   & 1.01    & 0.12                               \\
   & 4                     & $\SpectralLn$  & 0.95                          & \textbf{1.21} &     1.45 &      & 1.72      & 0.14                                                \\
   & 4                     & $\SpectralBn$ & 0.95                         & \textbf{1.21}      & \textbf{1.21}    &       & 2.06     & 0.15                                            \\\midrule
IEEE-118  & 4                     & Fastgreedy  & 1.11                            & \textbf{1.11}     & \textbf{1.11}   &   & 5.77   & 0.10                                  \\
  & 4                     & $\SpectralLn$  & 1.11                          & \textbf{1.14}     & 1.15  &      & 3.44           & 0.18                                        \\
  & 4                     & $\SpectralBn$ & 1.11                            & 1.16    & \textbf{1.11 }  &      & 2.91           & 0.17                                       \\\midrule
ACTIV-200 & 4                     & Fastgreedy  & 0.63                        & 0.72        & \textbf{0.69}   &         & 7.97    & 0.09                                           \\
 & 4                     & $\SpectralLn$  & 0.63                         & 0.72       & \textbf{0.63}   &             & 9.08      & 0.16                                     \\
 & 4                     & $\SpectralBn$ & 0.63                        & 0.71        & \textbf{0.63}   &          & 10.48           & 0.19                                 \\
\bottomrule
\end{tabular}
\vspace{-0.45cm} 
\end{table}

\section{Conclusion}
\label{sec:conclusion}
In this paper we considered an optimization problem to refine the bridge-block decomposition of power networks by means of line switching actions while having minimal impact on the maximum congestion. We revisited a heuristic two-stage tree partitioning procedure and proposed an MILP-based algorithm to solve its second stage using DC power flows. Numerical experiments on several test networks show that the improved two-stage approach performs better than the recursive approach, while overcoming the computational bottleneck of earlier brute-force algorithms. Furthermore, we proposed a modification of the recursive approach of~\cite{Zocca} to account for AC power flows. We show that, at least on five small instances, the recursive approach obtains objective values similar to the pre-existing brute-force two-stage approach, while drastically improving the running times.

Our numerical experiments suggest that the quality of the partitions plays an important role in achieving low post-switching congestion. As future work, we envision finding new optimization problem formulations to determine better partitions. Moreover, we plan to extend our current MILP for the OBS problem to account for AC power flows based on sophisticated linearization techniques in the literature \cite{Coffrin2014}.

\bibliographystyle{IEEEtran}
\bibliography{bibliography}
\end{document}

%% file: macros.tex
\usepackage{mathtools}
\usepackage{amsmath}
\usepackage{placeins}
\usepackage{cleveref}
\crefname{equation}{}{}
\usepackage{amssymb}
\usepackage{bm}
\DeclarePairedDelimiter{\abs}{\lvert}{\rvert}

\newcommand{\calE}{\mathcal{E}}

\newcommand{\calP}{\mathcal{P}}

\newcommand{\calV}{\mathcal{V}}

\newcommand{\bbR}{\mathbb{R}}

\newcommand{\defeq}{\vcentcolon=}

\newcommand{\nodeset}[1]{[#1]}

\usepackage{booktabs}
\usepackage{graphicx}
\usepackage{mathtools}

\newcommand{\crossE}{E_{C}(\calP)}

\newcommand{\binaryset}{\{0, 1\}}

\newcommand{\GcalE}{G^{\calE}}

\newcommand\clearrow{\global\let\rowmac\relax}
\clearrow

\newcommand{\SpectralBn}{\text{Spectral } B_N}
\newcommand{\SpectralLn}{\text{Spectral } L_N}
\usepackage{array}
\newcolumntype{H}{>{\setbox0=\hbox\bgroup}c<{\egroup}@{}}